\documentclass[a4paper,10pt]{amsart}
\usepackage[english]{certus}
\usepackage{a4wide}
\usepackage{amsmath,amsfonts,amssymb,amsthm}
\usepackage{comment}

\theoremstyle{definition}

\newtheorem*{As*}{Assumption}
\theoremstyle{plain}
\newtheorem*{Claim*}{Claim}
\theoremstyle{remark}
\DeclareMathOperator{\spa}{span}
\DeclareMathOperator{\sgn}{sgn}

\DeclareMathOperator{\conj}{conj}

\newcommand{\NN}{\ensuremath{\mathbb{N}}}
\newcommand{\ZZ}{\ensuremath{\mathbb{Z}}}
\newcommand{\RR}{\ensuremath{\mathbb{R}}}
\newcommand{\CC}{\ensuremath{\mathbb{C}}}

\title[On invariant random positive definite functions]{On invariant random positive definite functions}
\author{Vadim Alekseev}
\address{Vadim Alekseev, Technische Universit\"{a}t Dresden, Fachrichtung Mathematik, Institut f\"{u}r Geometrie, 01062 Dresden, Deutschland}
\email{vadim.alekseev@tu-dresden.de}

\author{Rahel Brugger}
\address{Rahel Brugger, Technische Universit\"{a}t Dresden, Fachrichtung Mathematik, Institut f\"{u}r Geometrie, 01062 Dresden, Deutschland}
\email{rahel.brugger@tu-dresden.de}

\subjclass[2010]{22D25, 20P99, 47A35, 46L55, 20C32}

\begin{document}
\onehalfspace

\begin{abstract}
	We give the definition of an invariant random positive definite function on a discrete group, generalizing both the notion of an invariant random subgroup and a character. We use von Neumann algebras to show that all invariant random positive definite functions on  groups with infinite conjugacy classes which integrate to the regular character are constant.
	
	This article is based on some results from the second named author's Ph.D. thesis.
\end{abstract}

\maketitle
\section{Introduction}

In the last years there has been a lot of progress about invariant random subgroups (IRSes), which shifted the attention in the study of ergodic group actions from their orbit equivalence relations to their stabilizers \cite{IRS12}, \cite{AGN2017}, \cite{Ge2018}, \cite{7s2011},\cite{7s2017}.
IRSes are a tool to study actions, but also behave similarly to normal subgroups.

We define a generalization of invariant random subgroups, which we call invariant random positive definite functions.
 An invariant random positive definite function (i.r.p.d.f.\@) is a measurable $\Gamma$-equivariant map
	\[\varphi\colon \Omega\to \mathrm{PD}(\Gamma),\]
	where $(\Omega,\mu)$ is a standard probabilitiy space with a measure preserving $\Gamma$-action, and $\mathrm{PD}(\Gamma)$ are the normalized positive definite functions $\phi$ on $\Gamma$ with $\Gamma$-action given by  $(g.\phi)(h)=\phi(g^{-1}hg)$ for $\phi\in\mathrm{PD}(\Gamma)$ and $g,h\in\Gamma$. This specializes to the definition of an IRS if we demand each $\varphi(\omega)$ to be the characteristic function of the stibilizer subgroup of $\omega$.

The definition of an i.r.p.d.f. is also closely related to the notion of a character on $\Gamma$, i.e.\@ a conjugation invariant normalized positive definite function. Indeed, if $\varphi$ is an i.r.p.d.f.\@, 
\[\mb E[\varphi]:=\int_{\Omega}\varphi(\omega)\,d\omega\]	
is a character.

A construction of Anatoly Vershik shows that in the case of $\Gamma=S_\infty$ every extremal character, except for the regular, the trivial and the alternating character, is of this form for a non-constant i.r.p.d.f.\@ $\varphi$ \cite{VerKer1981}. Some of these i.r.p.d.f.'s are IRSes, some are "twisted IRSes" arising from cocyles of the action.

Our main result is the following theorem. We call this phenomenon ``disintegration rigidity'' of the regular character $\delta_e\in\mathrm{Ch}(\Gamma)$.
\begin{Thm}[Theorem \ref{disintegration rigidity}]
	Let $\Gamma$ be a group where every nontrivial conjugacy class is infinite and let  $\varphi\colon \Omega\to \mathrm{PD}(\Gamma)$ be an i.r.p.d.f.~on $\Gamma$ with $\mb E[\varphi]=\delta_e$. Then $\varphi(\omega)=\delta_e$ for almost every $\omega\in\Omega$.
\end{Thm}

$\Gamma$ having infinite conjugacy classes is equivalent to $\delta_e\in\mathrm{Ch}(\Gamma)$ being an extremal character, hence the theorem states disintegration rigidity of $\delta_e$ in all cases where it has a chance to be disintegration rigid.

The main step in the proof of this theorem is to translate a given ergodic i.r.p.d.f.~$\varphi$ with $\mb E[\varphi]=\delta_e$ into a random variable $f\colon \Omega\to L^1(L\Gamma)$ which fulfills the invariance condition $f(\gamma.\omega)= \pi(\gamma^{-1})f(\omega)\pi(\gamma)$.
We then show that such a function must be constantly $1$, using that the conjugation action of $\Gamma$ on $L\Gamma$ is weakly mixing. Then $\varphi$ also must be constant.
This method might also apply to other characters than the regular one.

\subsection*{Acknowledgements}
The notion of an invariant random positive definite function was coined by Miklós Abért, the first named author and Andreas Thom at the Erwin Schrödinger Institute programme ``Measured group theory'', February 2016. Correspondingly, we would like to thank the Erwin Schrödinger Institute for the nice working atmosphere during the workshop; during it, both the authors were also partially supported by the ERC grant “ANALYTIC” no. 259527 of Goulnara Arzhantseva. We would also like to thank and Miklós Abért and Andreas Thom for sharing some nice ideas and discussing interesting questions around this project and Jesse Peterson for stimulating discussions and comments on an early version of this paper.

\section{Preliminaries}
As invariant random positive definite functions generalize both characters and invariant random subgroups we first collect some information about these.
\subsection{Characters on discrete groups}\label{section characters}
Let $\Gamma$ be a discrete, countable group.
\begin{Def}
	A function $\phi\colon\Gamma\to\CC$ is called \emph{positive definite} if for all $g_1,\dots,g_n\in\Gamma$ the matrix
	$ [\phi(g_j^{-1}g_i)]\in\mathbb{M}_n(\CC)$ is positive or, equivalently, if $\phi$ induces a state on $\CC\Gamma$.
\end{Def} 

\begin{Def}
	A \emph{character} $\tau\in\mathrm{Ch}(\Gamma)$ is a conjugation-invariant positive definite function on $\Gamma$  normalized by $\tau(e)=1$. A character is called \emph{extremal} if it is not a non-trivial convex combination of two different characters.
\end{Def}

The characters of a given group $\Gamma$ form a Choquet simplex, i.e.\@ every character can be uniquely decomposed as a convex combination of extremal ones \cite{Thoma1964(1)}.

If $\phi$ is a positive definite function, then $\left\langle g,h\right\rangle =\phi(h^*g)$ for $g,h\in\Gamma$ extends to a prescalar product on $\CC\Gamma$. Let $H$ be the separated completion and denote the image of $\delta_g$ in $H$ again by $\delta_g$. Then $\pi(g)\colon\delta_h\mapsto\delta_{gh}$ extends uniquely to a unitary operator $\pi(g)\in U(H)$. We get a unitary representation $\pi\colon\Gamma\to U(H)$ such that $\delta_e\in H$ is cyclic and
\[\phi(g)=\left\langle \pi(g)\delta_e,\delta_e\right\rangle \]
for all $g\in\Gamma$.
The triple $(H,\pi,\delta_e)$ is unique with these properties  up to a unitary. This is called the \emph{GNS construction} of $\phi$.
Sometimes we will also call the von Neumann algebra $\pi(\Gamma)''\subset B(H)$ the GNS construction of $\phi$.

If $\phi=\tau$ is a character, its GNS construction is a finite von Neumann algebra with trace extending the character. We denote this trace again by $\tau$ and get $L^2(\pi(\Gamma)'',\tau)=H$.
In this case we also have a unitary right representation 
\[\rho\colon \Gamma\to U(H), \quad \rho(g)\colon\delta_h\mapsto\delta_{hg^{-1}}.\]
Restricted to $\pi(\Gamma)''\subset L^2(\pi(\Gamma)'',\tau)$, the maps $\pi(g)$ and $\rho(g)$ correspond to $x\mapsto \pi(g)x$ and $x\mapsto x \pi(g^{-1})$ when $x$ is viewed as an operator $x\in B(H)$. In particular, 
\[\Gamma\to \Aut(\pi(\Gamma)''),\quad g\mapsto( x\mapsto \pi(g)x\pi(g^{-1})),\]
is a trace-preserving  action. 

In the case of the regular character $\delta_e$ we get the group von Neumann algebra $L\Gamma$ as GNS construction.

By \cite{Thoma1964(1)}, a character is extremal  if and only if its von Neumann algebra $\pi(\Gamma)''\subset  B(H)$ is a factor.

\begin{Def}
The type of a character is the type I of its von Neumann algebra (e.g. I, II$_1$ etc.). 
\end{Def}
Since the GNS construction of a character is finite, an extremal character can only be of type I$_n$ or II$_1$.
\subsection{Invariant random subgroups}\label{section IRS}
The name ``invariant random subgroup'' is due to \cite{IRS12}. However, the concept is much older and was, for example, studied by Vershik in the 80s and by Stuck-Zimmer in the 90s.
\begin{Def}\label{IRS}
	An invariant random subgroup (IRS) is a map given by
	\[\varphi\colon \Omega\to \mathrm{Sub}(\Gamma),\quad \omega\mapsto \mathrm{Stab}(\omega)=\{\gamma\in\Gamma\,|\,\gamma.\omega=\omega\},\]
	for a measure preserving action $\Gamma\curvearrowright (\Omega,\mu)$ on a standard probability space.
\end{Def}

In fact, invariant random subgroups were originally defined as conjugation invariant measures on $\mathrm{Sub}(\Gamma)$. One can show that this is equivalent to the above definition \cite[Proposition 13]{IRS12}.
We use this formulation because it will fit with our definition of invariant random positive definite functions and makes our notation easier.

If $\varphi\colon \Omega \to \mathrm{\mathrm{PD}}(\Gamma)$ is an IRS,
\[ \mathbb{ E}[\varphi]\colon \gamma\mapsto \mu(\{\omega|\,\gamma.\omega=\omega\}) \]
is a character.
\begin{Ex}\label{Sn}
	Let $\Omega=\{1,\dots,n\}$, let $\mu$ be the normalized counting measure and $\Gamma=S_n$ the symmetric group. Then $\varphi(i)=\{\sigma|\,\sigma(i)=i\}$ for $i\in\Omega$ is an IRS where $\mathbb{E}[\varphi]= \tr$ is the normalized trace on matrices. The trace $\tr$ is not an extremal character on $S_n$.
\end{Ex}
The following two theorems show that for $\Gamma=S_\infty$, many characters arise in this way.

\begin{Thm}[\cite{Thoma1964(3)}]\label{characters on Sinf}
	Every extremal character on $S_\infty$ is of the form
	\[\tau_{\alpha,\beta}(g)=\prod_{k\geq 2}s_k^{r_k(g)},	\]	
	where 	$r_k(g)$ is the number of cycles of length $k$ in $g$,
	$\alpha=(\alpha_n)_{n\in\NN}$ and $\beta=(\beta_n)_{n\in\NN}$ are sequences with 
	$\alpha_n\geq\alpha_{n+1}\geq 0$ and $\beta_n\geq\beta_{n+1}\geq 0$ for all $n\in\NN$ and such that  \[\sum_{n\in\NN}\alpha_n+\sum_{n\in\NN}\beta_n\leq 1\] and the $s_k$ are given by
	\[s_k:=\sum_{n\in\NN}\alpha_n^k+(-1)^{k+1}\sum_{n\in\NN}\beta_n^k.\]
	All such $\tau_{\alpha,\beta}$ are extremal characters and $\tau_{\alpha,\beta}=\tau_{\alpha',\beta'}$ implies $\alpha=\alpha'$ and $\beta=\beta'$.
	
	All extremal characters on $S_\infty$ exept for the trivial character and the alternating character are of type II. 
\end{Thm}
\begin{Rem}\label{trivial characters on Sinf}
	In the theorem the trivial character belongs to $\alpha=(1,0,0,\dots)$ and $\beta=0$, the alternating character belongs to $\alpha=0$ and $\beta=(1,0,0,\dots)$ and the regular character belongs to $\alpha=\beta=0$.
\end{Rem}
\begin{Thm}[\cite{VerKer1981}]\label{tau on S as IRS}
	Using the notation of Theorem \textup{\ref{characters on Sinf}}, assume $\beta=0$, let  \[\delta=1-\sum_{n\in\NN}\alpha_n\] and
	let $\mathcal{Q}=\NN\sqcup [0,\delta]$ with probability measure $\mu$ which is $(\alpha_n)_{n\in\NN}$ on $\NN$ and the Lebesgue measure on $[0,\delta]$.
	Let $\Omega=\prod_{1}^\infty\mathcal{Q}$ with measure $m_{\alpha,0}=\prod_1^\infty\mu$ and let $S_\infty$ act on $(\Omega, m_{\alpha,0})$ by permutation of the coordinates.
	
	Then $\tau_{\alpha,0}=\mathbb{E}[\varphi]$ for this IRS $\varphi$.
\end{Thm}

\section[Definition and examples of i.r.p.d.f.'s]{Invariant random positive definite functions}
\label{section irpdf}
\begin{Def}\label{irpdf}
	Let $\Gamma$ be a discrete group. An invariant random positive definite function (i.r.p.d.f.\@) is a measurable $\Gamma$-equivariant map
	\[\varphi\colon \Omega\to \mathrm{PD}(\Gamma),\]
	where $(\Omega,\mu)$ is a standard probabilitiy space with a measure preserving $\Gamma$-action and $\mathrm{PD}(\Gamma)$ are the positive definite functions $\phi$ on $\Gamma$ with $\phi(e)=1$ and $\Gamma$-action given by $(g.\phi)(h)=\phi(g^{-1}hg)$ for $\phi\in\mathrm{PD}(\Gamma)$ .
\end{Def}
We often write $\varphi_\omega$ for $\varphi(\omega)$. 
\begin{Def}
	An i.r.p.d.f.~$\varphi$ is called \emph{ergodic} if the action $\Gamma\curvearrowright (\Omega,\mu)$ is ergodic. 
	
	We say $\varphi$ is \emph{extremal} if $\varphi=c\varphi_1+(1-c)\varphi_2$ for  i.r.p.d.f.'s $\varphi_i\colon\Omega\to\mathrm{PD}(\Gamma)$ and $c\in (0,1)$ implies that $\varphi_1=\varphi_2=\varphi$.
\end{Def}

When viewing the i.r.p.d.f.'s with given $\Gamma\curvearrowright\Omega$ as $\Gamma$-equivariant positive definite functions $\varphi\colon \Gamma\to L^\infty(\Omega,\mu)$, they form a compact convex subset of $\ell^\infty(\Gamma,L^\infty(\Omega,\mu))$ with the topology of pointwise weak$^*$ convergence.  By the Kre\u{\i}n--Milman Theorem, the space of these functions is then equal to the closed convex hull of its extremal points. Hence as for characters, every i.r.p.d.f.~is the convex integral of extremal i.r.p.d.f.'s.
\begin{Ex}
	Invariant random subgroups are i.r.p.d.f.'s because the subgroups $\mathrm{Sub}(\Gamma)$ of $\Gamma$ are canonically embedded in $\mathrm{PD}(\Gamma)$ by taking the characteristic function and the stabilizers of an action fulfill the invariance condition in Definition \ref{irpdf}.
\end{Ex}

As for invariant random subgroups, if $\varphi\colon \Omega \to \mathrm{\mathrm{PD}}(\Gamma)$ is an i.r.p.d.f.\@, 
\[ \mathbb{E}[\varphi] = \int_{\Omega} \varphi_\omega\,d\mu(\omega)\]
is a character.

\begin{Qu}
	Does ergodicity and extremality of $\varphi$ imply that $\mathbb{E}[\varphi]$ is extremal as a character?
\end{Qu}
A positive answer to this question would mean that it every i.r.p.d.f.\@ can be decomposed into i.r.p.d.f.'s  with an extremal character as expectation.
\begin{Ex}\label{sphere}
	Let $(S,\lambda)$ be the unit sphere in $\CC^n$ with Lebesgue measure and let $\Gamma$ be a discrete subgroup of the unitary group $U(n)$ acting on $S$ in the natural way. Then
	\[\varphi\colon S\to \mathrm{PD}(\Gamma),\quad \varphi_\xi(\gamma)=\left\langle \gamma.\xi,\xi\right\rangle \quad \forall\xi\in S,\gamma\in\Gamma	\]
	is an i.r.p.d.f.~for which $\mathbb{E}[\varphi]=\tr$ is the normalized trace on matrices, which is an extremal character on $\Gamma$ iff $\Gamma$ generates $\mb M_n(\CC)$ as an algebra. For such $\Gamma$, $\varphi$ is an extremal i.r.p.d.f..
\end{Ex}
\begin{Ex}\label{nonextremal character}
	Let $(S^1,\lambda)$ be the circle with Lebesgue measure and trivial action of $\ZZ$. Then
	\[\varphi\colon S^1\to \mathrm{PD}(\ZZ),\quad\varphi_z(n)=z^n\]
	is an i.r.p.d.f.~with $\mathbb{E}[\varphi]=\delta_e$. Here $\delta_e$ is not extremal and $\varphi_z$ is an extremal character for every $z\in S^1$. In this way every decomposition of a non-extremal character into extremal ones gives an i.r.p.d.f.~with trivial action.
\end{Ex}
\begin{Ex}\label{many phis}
	Let $G$ be a compact group with Haar measure $\mu$ and $\Gamma < G$. Let $\Gamma$ act on $G$ by left multiplication. Let $\pi\colon G\to U(H)$ be  a unitary representation and $\xi\in H$ a unit vector. Then 
	\[\varphi^\xi\colon (G,\mu)\to \mathrm{PD}(\Gamma),\quad \varphi^\xi_g(h)=\left\langle \pi(hg)\xi, \pi(g)\xi\right\rangle  	\]
	is an i.r.p.d.f..	
	If $\pi\colon G\to U(\CC^n)$ is irreducible and $\Gamma$ is dense, then $\mathbb{E}[\varphi^\xi](\gamma)=\tr(\pi(\gamma))$, which is an extremal character on $\Gamma$, and $\varphi^\xi$ is ergodic and extremal.
\end{Ex}
Example \ref{many phis} shows that, in contrast to the situation for characters, the decomposition of an i.r.p.d.f.~into extremal i.r.p.d.f.'s is not unique: Take an irreducible representation $\pi\colon G\to U(\CC^n)$, an orthonormal basis $(\xi_i)$ of $\CC^n$ and $\Gamma<G$ dense. Then
\[\sum_{i=1}^n \frac{1}{n}\varphi^{\xi_i}\equiv\tr\circ\pi.\]
For different bases we get different $\varphi^{\xi_i}$'s, so this gives different convex decompositions of the constant i.r.p.d.f.~$\tr\circ\pi$ into extremal i.r.p.d.f.'s..

\begin{Thm}[\cite{VerKer1981}, Theorem 3]\label{tau on S as irpdf}
	In the notation of Theorem \textup{\ref{characters on Sinf}}, let 
	\[\delta=1-\sum_{n\in\NN}\alpha_n-\sum_{n\in\NN}\beta_n,\]
	$\NN_+=\NN_-=\NN$ and $\mathcal{Q}=\NN_+\sqcup \NN_-\sqcup[0,\delta]$ with the probability measure $\mu$ which is $(\alpha_n)_{n\in\NN}$ on $\NN_+$, $(\beta_n)_{n\in\NN}$ on $\NN_-$ and the Lebesgue measure on $[0,\delta]$.
	Then let $\Omega=\prod_{1}^\infty\mathcal{Q}$ with the measure $m_{\alpha,\beta}=\prod_1^\infty\mu$ and let $S_\infty$ act on $(\Omega,m_{\alpha,\beta})$ by permutation of the coordinates.
	
	For $g\in S_\infty$ and $\omega\in\Omega$ define $\sgn(g,\omega)$ to be $1$ if 
	\[\prod_{(i,j) :\, \omega_i,\omega_j\in \NN_-,\, i<j}\left( g(j)-g(i)\right) 	\]
	is positive and $-1$ otherwise. This fulfills the cocycle identity
	\begin{align}\label{cocycle identity}
		\sgn(gh,\omega)=\sgn(h,\omega)\sgn(g,h.\omega).
	\end{align}
	Let \[\varphi_\omega(g)
	=\begin{cases} \sgn(g,\omega)\quad & \text{ if } g .\omega=\omega,\\
	0 \quad &\text{ if } g .\omega\ne\omega.
	\end{cases}
	\]	
	Then $\tau_{\alpha,\beta}=\mathbb{E}[\varphi]$.
\end{Thm}
The following theorem proves that the above $\varphi$ is an i.r.p.d.f.. If $\beta$ is non-trivial, then $\varphi$ is not an IRS.
\begin{Thm}\label{irpdf by c}
	Let $\Gamma\curvearrowright(\Omega,\mu)$ be a p.m.p.\@ action and $c\colon \Gamma\times\Omega\to S^1$ a cocycle as in \eqref{cocycle identity}. Then 
	\[\varphi_\omega(g)
	=\begin{cases} c(g,\omega)\quad & \text{ if } g .\omega=\omega,\\
	0 \quad &\text{ if } g .\omega\ne\omega,
	\end{cases}
	\]	
	is an i.r.p.d.f..
\end{Thm}
If $c$ is not constantly $1$, $\varphi$ is not an IRS because it takes values outside $\{0,1\}$.
\begin{proof}	
	To show that $\varphi$ is invariant we need that $c(g,h\omega)=c(h^{-1}gh,\omega)$ if $h^{-1}gh.\omega=\omega$. By the cocycle identity we have
	\[1=c(1,\omega)=c(h^{-1}h,h^{-1}gh.\omega)=c(h,\omega)c(h^{-1},gh.\omega)	\]
	and hence 
	\[c(h^{-1}gh,\omega)=c(h,\omega)c(h^{-1}g,h.\omega)=c(h,\omega)c(g,h.\omega)c(h^{-1},gh.\omega)=c(g,h.\omega).	\]
	
	Now we show that $\varphi_{\omega}$ is positive definite for a.e.\@ $\omega\in\Omega$.
	Let $\mathcal{R}\subset\Omega\times\Omega$ be the orbit equivalence relation of $\Gamma\curvearrowright(\Omega,\mu)$, equipped with the measure $\mu_{\mathcal{R}}$ which is $\mu$ on $\Omega$ and the counting measure in each fiber, i.e., for 	$A\subset \mathcal{R}$ measurable 
	\[\mu_{\mathcal{R}}(A):=\int_\Omega |\{(x,y)\in A\}| dx.	\]  
	 Then 
	$\pi\colon \Gamma \to U(L^2(\mathcal{R})),$ given by
	\[(\pi(g)\xi)(x,y)=c(g,x)\xi(g.x,y)	 \]
	is a unitary representation and for every $X\subset\Omega$ we find a vector $\xi_X=\chi_{\{(x,x)|x\in X\}}\in L^2(\mathcal{R})$ such that
	\[\int_X \varphi_{\omega}(g)=\left\langle\pi(g)\xi_X,\xi_X\right\rangle .\]
	Hence for every $a\in\CC\Gamma$ we have
	\[\int_X \varphi_{\omega}(a^*a)\geq 0\]
	for all $X\subset\Omega$ and hence $\varphi_{\omega}(a^*a)\geq 0$ almost everywhere.
\end{proof}
Up to now, all our examples of i.r.p.d.f.'s which integrate to a type II character are of this form. In particular, they
 are supported on an IRS in the sense that $\varphi_{\omega}(\gamma)=0$ if $\gamma.\omega\ne\omega$. This leads to the following questions.
\begin{Qu}
	Is every i.r.p.d.f.\@ $\varphi$ such that $\mathbb{E}[\varphi]$ is of type II$_1$  supported on an IRS ?
\end{Qu}
\begin{Qu}
	Is every i.r.p.d.f.\@ which is supported on an IRS  as in Theorem \ref{irpdf by c} ?
\end{Qu}
\section[I.r.p.d.f.'s and von Neumann algebras]{Connections to von Neumann algebras}
\label{general irpdf}
In this section we translate i.r.p.d.f.'s into the language of von Neumann algebras in order to be able to use von Neumann methods to study them in the next section.
For the relevant theory of von Neumann algebras see \cite{Blackadar2006},\cite{ADPopa},\cite{Houdayer}.

Fix a discrete group $\Gamma$,  a character $\tau\in \mathrm{Ch}(\Gamma)$ and  an ergodic, measure preserving action $\alpha\colon\Gamma\curvearrowright (\Omega,\mu)$ on a standard probability space. Let $A:=L^\infty(\Omega,\mu)$ and write again $\alpha$ for the corresponding action on $A$.
Let $\pi\colon \Gamma\to U(H)$ be the GNS representation of $\tau$. 

\begin{Lemma}\label{M varphi}
	Let $\varphi$ be an i.r.p.d.f.~with $\mb E[\varphi]=\tau$ and for each  $\omega\in\Omega$ let $(\pi_{\omega},H_\omega,\xi_\omega)$ be the GNS construction of $\varphi_\omega$. Let 
	\[H_\varphi:=\int_\Omega^\oplus H_\omega\, d\mu(\omega)\]
	be the direct integral of Hilbert spaces, $\xi=(\xi_\omega)_{\omega\in\Omega}\in H_\varphi$ and \[\pi_\varphi=\int_\Omega^\oplus\pi_\omega\, d\mu(\omega)\colon \Gamma\to  B(H_\varphi)\] the direct integral of representations.
	Then $\pi_{\varphi}(\Gamma)''\cong\pi(\Gamma)''$ with isomorphism taking $\pi_\varphi(\gamma)$ to $\pi(\gamma)$ for all $\gamma\in\Gamma$.
\end{Lemma}
\begin{proof}
	Let $p\in  B(H_\varphi)$ be the orthogonal projection onto the cyclic representation of $\xi$. Then $p\in\pi_\varphi(\Gamma)'$.
	As $\mb E[\varphi]=\tau$ we have 
	\[\left\langle \pi_\varphi(\gamma)\xi,\xi\right\rangle 
	=\int_\Omega \left\langle \pi_\omega(\gamma)\xi_\omega,\xi_\omega\right\rangle \,d\mu(\omega)
	= \int_\Omega \varphi_\omega(\gamma)\,d\mu(\omega)=\tau(\gamma) \]
	for all $\gamma\in\Gamma$.
	So $(p(H_\varphi),\pi_\varphi,\xi)$ is a GNS triple for $\tau$ and therefore by uniqueness of the GNS construction \[\pi(\Gamma)'' \cong (p\,\pi_\varphi(\Gamma)\,p)''=p (\pi_\varphi(\Gamma))''\]
	with isomorphism taking $\pi(\gamma)$ to $p\pi_\varphi(\gamma)$ for all $\gamma\in\Gamma$.
	Now we show that \[\Phi\colon (\pi_\varphi(\Gamma))''\to p (\pi_\varphi(\Gamma))'',\; x\mapsto px,\]is an isomorphism. It is clearly a surjective homomorphism. For injectivity let $x\in (\pi_\varphi(\Gamma))''$ with $\Phi(x^*x)=px^*xp=0$. Then for all $a\in\CC\Gamma$ we have
	\[0=\left\langle x^*x \pi_\varphi(a)\xi,\pi_\varphi(a)\xi\right\rangle 
	=\int_\Omega \left\langle (x^*x)_\omega \pi_\omega(a)\xi_\omega,\pi_\omega(a)\xi_\omega\right\rangle d\mu(\omega)
	\]
	and therefore $\left\langle (x^*x)_\omega \pi_\omega(a)\xi_\omega,\pi_\omega(a)\xi_\omega\right\rangle=0$ for almost all $\omega\in\Omega$. But  $\pi_\omega(\CC\Gamma)\xi_\omega$ is dense in $H_\omega$, so $(x^*x)_\omega=0$ for almost all $\omega$. Hence $x=0$ and $\Phi$ is injective.
	
	Composing the two isomorphisms we get $\pi(\Gamma)''\cong p\pi_\varphi(\Gamma)''\cong \pi_\varphi(\Gamma)''$ with isomorphism mapping $\pi(\gamma)$ to $\pi_{\varphi}(\gamma)$.
\end{proof}

\begin{Lemma}
	Let $M := (A\cup \pi_\varphi(\Gamma))''=\int_\Omega^\oplus\pi_\omega(\Gamma)''d\mu(\omega)$.	Then $M$ is a finite von Neumann algebra.
\end{Lemma}
\begin{proof}
	Let $u\in M$ be such that $u^* u = 1$. By the Kaplansky Density Theorem we find a sequence of finite sums
	\[
	t_n = \sum_{i} p_{n,i} x_{n,i}
	\]
	converging to $u$ in the strong$^\ast$ topology such that $\|t_n\|\leq1$ for all $n$, $p_{n,i}\in A$ are mutually orthogonal projections  for fixed $n$ and $x_{n,i}\in \pi_\varphi(\Gamma)''$. We then have $t_n^*t_n \overset{s^*}{\to} 1$  since the strong* topology is jointly continuous on bounded sets. Hence  $|t_n| \overset{s^*}{\to} 1$ by \cite[Lemma II.4.6]{ Takesaki2002}. Letting
	\[
	f(t)\coloneqq
	\begin{cases}
	1-2t,& 0\leqslant t\leqslant 1/2,\\
	0, & 1/2\leqslant t \leqslant 1,
	\end{cases}
	\]
	we obtain (again by \cite[Lemma II.4.6]{ Takesaki2002}) $f(|t_n|) \overset{s^*}{\to} 0$, and therefore $|t_n| + f(|t_n|) \overset{s^*}{\to} 1$. However, as $1/2\leqslant t + f(t)\leqslant 1$ on $[0,1]$, we also have $1/2 \leqslant |t_n| + f(|t_n|) \leqslant 1$.
	
	Let $t_n = u_n |t_n|$ be the polar decomposition of $t_n$. Then we have
	\[
	u_n (|t_n| + f(|t_n|)) \overset{s^*}{\to} u
	\]
	because $f(|t_n|) \overset{s^*}{\to} 0$. On the other hand, $|t_n| + f(|t_n|)$ is invertible with the inverse bounded by $2$ and $(|t_n| + f(|t_n|))^{-1} \overset{s^*}{\to} 1$ again by \cite[Lemma II.4.6]{ Takesaki2002}. Therefore,
	\begin{align}\label{u_n to u}
	u_n = u_n (|t_n| + f(|t_n|)) (|t_n| + f(|t_n|))^{-1} \overset{s^*}{\to} u.
	\end{align}
	
	Let $x_{n,i}=v_{n,i}|x_{n,i}|$ be the polar decomposition of $x_{n,i}$. Then
	\begin{align}\label{u_n}
	u_n = \sum_{i} p_{n,i} v_{n,i}
	\end{align}
	because using that $A$ commutes with $\pi_\varphi(\Gamma)''$ and that the $p_{n,i}$ are mutually orthogonal we get that
	\[|t_n|=\sum_ip_{n,i}|x_{n,i}|,	\]
	and hence 
	\[\left(  \sum_{i} p_{n,i} v_{n,i}\right) |t_n|
	= \sum_{i} p_{n,i} v_{n,i}|x_{n,i}|
	= \sum_{i} p_{n,i}x_{n,i}=t_n.	\]
	
	Now (\ref{u_n}) and (\ref{u_n to u}) imply that
	\[
	u_n^*u_n = \sum_{i} p_{n,i} v_{n,i}^*v_{n,i} \overset{s^*}{\to} 1,
	\]
	and therefore
	\begin{align}\label{sum p}
	\sum_{i} p_{n,i} \overset{s^*}{\to} 1.
	\end{align}
	
	Since $\pi_\varphi(\Gamma)''\cong\pi(\Gamma)''$ is finite, there exist partial isometries $w_{n,i}\in\pi_\varphi(\Gamma)''$ such that $u_{n,i} = v_{n,i} + w_{n,i}$ are unitaries. Let $q_{n,i}\coloneqq w_{n,i}^*w_{n,i}$ be the source projections of the $ w_{n,i}$. Then
	\[
	\sum_{i} p_{n,i} q_{n,i}
	= \sum_{i} p_{n,i}(1-v_{n,i}^*v_{n,i})
	\leq 1-\sum_{i} p_{n,i}v_{n,i}^*v_{n,i}
	=1 - u_n^*u_n \overset{s^*}{\to} 0,
	\]
	and therefore
	\[
	\sum_{i} p_{n,i} w_{n,i} = \left(\sum_{i} p_{n,i} w_{n,i}\right)\left(\sum_{i} p_{n,i} q_{n,i}\right) \overset{s^*}{\to} 0.
	\]
	Thus by (\ref{u_n})
	\[
	\sum_{i} p_{n,i} u_{n,i} = u_n + \sum_{i} p_{n,i} w_{n,i} \overset{s^*}{\to} u,
	\]
	and therefore, since the $u_{n,i}$ are unitaries,
	\[
	\sum_{i} p_{n,i} =\sum_{i} p_{n,i} u_{n,i}u_{n,i}^* \overset{s^*}{\to} uu^*.
	\]
	Hence $uu^* = 1$ by (\ref{sum p}), which means that $M$ is finite.
\end{proof}

\begin{Lemma}\label{M tensor product}
	If $\tau$ is extremal, we have $M \cong A\vntens \pi_\varphi(\Gamma)''$ with isomorphism taking $xa\in M$ to $a\otimes x\in A\vntens \pi_\varphi(\Gamma)''$ for all $a\in A$ and $x\in\pi_\varphi(\Gamma)''$.
\end{Lemma}
\begin{proof}
	Since $M$ is finite by the previous lemma, there exists a normal faithful conditional expectation $E\colon M\to \pi_\varphi(\Gamma)''$. Since $\pi_\varphi(\Gamma)''$ and $A$ commute and $E$ is $\pi_\varphi(\Gamma)''$-linear, 
	\[
	E(a) = E(\pi_\varphi(\gamma) a \pi_\varphi(\gamma^{-1})) = \pi_\varphi(\gamma) E(a) \pi_\varphi(\gamma^{-1})
	\]
	for all $ \gamma\in \Gamma$ and $ a\in A.$ Thus, $E(A)$ is contained in the center of $\pi_\varphi(\Gamma)''\cong\pi(\Gamma)''$, which is equal to $\CC$ since $\tau$ is extremal. Now the claim follows from \cite[Theorem 9.12]{Stratila1981}.
\end{proof}
On $M$  resp. $L^1(M)$ we define a $\Gamma$-action $\theta$ by 
\[\theta_\gamma(a\otimes m)=\alpha_\gamma(a)\otimes \pi(\gamma)m\pi(\gamma^{-1}).\]
By  $M^\theta$ resp. $ L^1(M)^\theta$ we denote the elements that are invariant under $\theta$.
\begin{Prop}\label{M theta}
	Given an ergodic action and an extremal character $\tau\in\mathrm{Ch}(\Gamma)$ there is a one-to-one correspondence between i.r.p.d.f.'s $\varphi\colon \Omega\to \mathrm{PD}(\Gamma)$ with $\mb E[\varphi]=\tau$ and positive selfadjoint elements $ f\in L^1(M)^\theta$ with $\int_{\Omega}f_\omega\, d\mu(\omega)=1$ such that
	\[\varphi_\omega(\gamma)	=\tau(\pi(\gamma)f_\omega).\]
\end{Prop}
\begin{proof}
	By Lemma \ref{M varphi} and Lemma \ref{M tensor product}, we have $\pi(\Gamma)''\cong \pi_\omega(\Gamma)''$ for a.e.\@ $\omega\in\Omega$ with the canonical isomorphism  sending $\pi(\gamma)$ to $\pi_\omega(\gamma)$ for each $\gamma\in\Gamma$.
	As $\varphi_\omega(\gamma)=\left\langle \pi_\omega(\gamma)\xi_\omega, \xi_\omega\right\rangle $, we can extend it to
	\[\varphi_\omega\colon \pi_\omega(\Gamma)''\to \CC,\;x\mapsto \left\langle x\xi_\omega,\xi_\omega\right\rangle, \]
	which is a positive normal functional on $\pi_\omega(\Gamma)''$ and therefore on $\pi(\Gamma)''$. So by  \cite[Lemma IX.2.12]{Takesaki} there exists a unique positive element $f_\omega\in L^1(\pi(\Gamma)'',\tau)$ such that $\varphi_\omega(x)=\tau(xf_\omega)$ for all $x\in\pi(\Gamma)''$. 	
	Let $f\colon\Omega\to L^1(\pi(\Gamma)''),\,\omega\mapsto f_\omega$.
	To see that $f$ is $\theta$-invariant, we calculate
	\[\tau(\pi(\gamma)f_{\alpha_{\gamma'}(\omega)})
	=\varphi_{\alpha_{\gamma'} (\omega)}(\gamma)
	=\varphi({\gamma'}^{-1} \gamma\gamma')
	=\tau(\pi(\gamma)\pi(\gamma')f_{\omega} \pi({\gamma'}^{-1})),	\]
	so $\alpha^{-1}_{\gamma'}(f)_\omega=f_{\alpha_{\gamma'}(\omega)}=\pi(\gamma')f_\omega\pi({\gamma'}^{-1})$ for all $\gamma'\in \Gamma$ by uniqueness of $f$, hence $\theta(f)=f$.
	It follows that $\|f_\omega\|_1$ is $\Gamma$-invariant and hence constant, so $f\in L^1(M)^\theta$.
	We have for all $\gamma\in\Gamma$
	\[ \tau\left( \pi(\gamma) \int f_\omega\, d\mu(\omega)\right) 
	=\int \tau\left( \pi(\gamma) f_\omega \right)\, d\mu(\omega)
	=\int \varphi_\omega(\gamma)\,d\mu(\omega)
	=\tau(\pi(\gamma)),	\]
	hence $\int f_\omega\, d\mu(\omega)=1$.
	By \cite[Lemma 8.3 (3)]{Luck}, $f$ is a selfadjoint operator. 
	
	Conversely it is easy to check that such an $f$ defines an i.r.p.d.f.~$\varphi$ with $\mb E(\varphi)=\tau$ by $\varphi_\omega(\gamma)=\tau(\pi(\gamma)f_\omega)$.
\end{proof}
\begin{Rem}\quad
	\begin{enumerate}
		\item 	If $\varphi$ is as in Example \ref{sphere} with $\Gamma$ big enough so that $\varphi$ is extremal, we have $f\colon S\to \mathbb{M}_n(\CC)$ with $f_\xi$ the orthogonal projection on $\spa (\xi)$.
		\item 	Similarly, if $\Gamma$ in Example \ref{many phis} is dense and $\pi$ irreducible, we find $f\colon G\to\mathbb{M}_n(\CC)$ where $f_g$ is the orthogonal projection on $\spa (\pi(g)\xi)$.
		\item 	The i.r.p.d.f.~in Example \ref{nonextremal character} is not of the form as in Proposition \ref{M theta}. Hence the ergodicity and extremality assumptions are necessary (or at least one of them is).
	\end{enumerate}
\end{Rem}
\begin{Lemma}
	In fact, for $f\in L^1(M)^\theta$ as in Proposition \textup{\ref{M theta}} the condition that $\int_{\Omega}f_\omega\, d\mu(\omega)=1$ is equivalent to 
	$\tau_M(f)=1$, where $\tau_M=\int_{\Omega}\otimes\,\tau$ is the trace on $M$.
\end{Lemma}
\begin{proof}
	Let $f$ be constructed from $\varphi$ as above.	Then
	\[\tau_M(f)=\int_{\Omega}\tau(f_\omega)\,d\mu(\omega)=\int_{\Omega}\varphi_\omega(e)\,d\mu(\omega)=\int_{\Omega}1\,d\mu(\omega)=1.	\]
	For the other direction let first $p\in M^\theta$ be a projection. Then 
	\[\tau(\gamma)=\tau\left( \pi(\gamma) \int p_\omega \,d\mu(\omega)\right) +\tau\left( \pi(\gamma) \int(1-p)_\omega\, d\mu(\omega)\right)  \] 
	is a convex decomposition into two characters. So by extremality of $\tau$, 
	\[\int p_\omega \,d\mu(\omega)=\tau_{M}(p)\cdot 1.\]
	Now let $f\in L^1(M)^\theta$ be positive selfadjoint with $\tau_M(f)=1$. Then it follows from the above and the spectral theorem for $f$ that  
	$\int f_\omega\, d\mu(\omega)=\tau_{M}(f)\cdot 1=1$.
\end{proof}
\begin{Lemma}	For $\tau$ extremal  and $\alpha$ ergodic
	the extremal i.r.p.d.f.'s $\varphi$ given $\alpha$ and $\mb E[\varphi]=\tau$ correspond to minimal projections in $M^\theta$. 
	$M^\theta$ is a direct sum of matrix algebras.
\end{Lemma}
\begin{proof}
	Let $\varphi\colon \Omega\to \mathrm{PD}(\Gamma)$ be an extremal i.r.p.d.f.~and $f\in L^1(M)^\theta$ as in Proposition \ref{M theta} such that $\tau(f_\omega \pi(\gamma))=\varphi_\omega(\gamma)$ for a.e.\@ $\omega\in\Omega$ and all $\gamma\in\Gamma$. Assume that $f$ is not a scalar multiple of a projection. Then there is a $c\in\RR^+$ such that
	\[f^{<c}:=\chi([0,c))f \quad \text{and}\quad f^{\geq c}:=\chi([c,\infty))f,\]
	are both nonzero with $\chi(I)$ denoting the spectral projection on $I$. These are again positive elements in $M^\theta$ hence $\tau_M(f^{<c})^{-1}f^{<c}$ and $\tau_M(f^{\geq c})^{-1}f^{\geq c}$ define two different i.r.p.d.f.'s $\varphi^{<c}$ and $\varphi^{\geq c}$ such that
	\[\varphi=\tau_M(f^{<c})\varphi^{<c}+\tau_M(f^{\geq c})\varphi^{\geq c}\]
	contradicting the extremality of $\varphi$. So $f=\tau_M(p)^{-1}p$ for some projection $p\in M^\theta$.
	If $p$ is not minimal in $M^\theta$, say $q<p$ and $q\in M^\theta$, then again $q$ and $p-q$ define two i.r.p.d.f.'s such that a convex combination gives $\varphi$, which contradicts extremality.
	
	Conversely every minimal projection $p\in M^\theta$ gives an extremal i.r.p.d.f.~$\varphi$ because if there was a decomposition $\varphi=c\varphi_1+(1-c)\varphi_2$ for some $0<c<1$ and different i.r.p.d.f.'s $\varphi_i$, this would give different positive elements $f_1,f_2\in M^\theta$ such that 
	$\tau_M(p)^{-1}p=cf_1+(1-c)f_2$, which is not possible for a minimal projection $p$.
	
	Since the set of i.r.p.d.f.'s is the closed convex hull of its extremal points, every positive trace 1 element of $M^\theta$ is a convex integral of minimal projections. This means $M^\theta$ is generated by its minimal projections, hence it is of type I with no diffuse part, i.e., $\mathcal{Z}(M^\theta)=L^\infty(X,\mu)$ such that every point in $X$ has positive mass. Since it is also finite, it follows that $M^\theta$ is a (maybe infinite) direct sum of matrix algebras.
\end{proof}
\begin{Rem}
	Let $\tau\in\mathrm{Ch}(\Gamma)$ be an extremal character, $\alpha\colon\Gamma\curvearrowright \Omega$ an ergodic action and $\theta$ corresponding to $\alpha$ and $\tau$ as in Proposition \ref{M theta}. Then, for i.r.p.d.f.'s associated to $\alpha$, we have the following observations.
	\begin{enumerate}
		\item As $M^\theta$ is a direct sum of matrix algebras every i.r.p.d.f.\@ $\varphi$ with $\mb E[\varphi]=\tau$ is a convex combination of countably many extremal ones.
		\item 	$M^\theta=\CC$ iff the constant i.r.p.d.f.\@ $\tau$ is the only one  with $\mb E[\varphi]=\tau$. It is also equivalent to the constant $\tau$ being an extremal i.r.p.d.f..
		If this is true for all $\alpha$, $\tau$ is disintegration rigid.
		\item	$M^\theta$ is abelian iff the decomposition of i.r.p.d.f.'s with $\mb E[\varphi]=\tau$ into extremal ones is unique. 
		\item $M^\theta$ is finite-dimensional iff every i.r.p.d.f.~is a finite convex sum of extremal ones.
	\end{enumerate}
\end{Rem}

\section{Disintegration rigidity of the regular character on i.c.c.\@ groups}\label{section icc}
In this section we show the following theorem.
\begin{Thm}\label{disintegration rigidity}
	Let $\Gamma$ be a group with infinite conjugacy classes. Let $\varphi\colon \Omega\to \mathrm{PD}(\Gamma)$ be an i.r.p.d.f.\@ on $\Gamma$ with $\mb E[\varphi]=\delta_e$. Then $\varphi(\omega)=\delta_e$ for almost every $\omega\in\Omega$.
\end{Thm}
\begin{Def}
	If the conclusion of the theorem  holds, we say  $(\Gamma,\delta_e)$ is \textit{disintegration rigid}.
\end{Def}
\begin{Rem}
	Theorem \ref{characters on Sinf}, Remark \ref{trivial characters on Sinf} and Theorem \ref{tau on S as irpdf} show that the regular character, the trivial character and the alternating character are the only disintegration rigid characters on $S_\infty$. Indeed, if $S_\infty\curvearrowright (\Omega,m_{\alpha,\beta})$ is the action from Theorem \ref{tau on S as irpdf} such that $\tau_{\alpha,\beta}$ is none of these three characters, we have $0<\alpha_1<1$ or $0<\beta_1<1$. Assume w.l.o.g. that $0<\alpha_1<1$. Then for every nontrivial $g\in S_\infty$ and $j\in\supp(g)=\{j|\,g(j)\ne j\}$
	\begin{align*}
	0&<	m_{\alpha,\beta}\left( \{\omega\in\Omega|\,\omega_i=1\in\NN_+\;\forall i\in\supp(g)\}\right) \\
	&\leq m_{\alpha,\beta}\left( \{\omega\in\Omega|\,g.\omega=\omega\}\right) \\
	&\leq 1- m_{\alpha,\beta}\left( \{\omega\in\Omega|\,\omega_j=1\in\NN_+,\; g.\omega_j\ne 1\in\NN_+\}\right)
	<1.
	\end{align*}
	Hence the $\varphi$ in Theorem \ref{tau on S as IRS} is non-constant with $\mb E[\varphi]=\tau_{\alpha,\beta}$.
	
	The trivial and the alternating character are clearly disintegration rigid because every positive definite function takes values in the unit disk, and thus, if an i.r.p.d.f.~intergrates to a character which takes values only on the boundary of the unit disk, the i.r.p.d.f.~has to be constant.
\end{Rem}

\begin{Def}
	A trace-preserving action on a finite von Neumann algebra $\Gamma\to \Aut(M)$ is called \textit{weakly mixing} if $\CC\cdot 1$ is the only finite-dimensional, $\Gamma$-invariant subspace in $M$.
\end{Def}
The following lemma might be known to experts but we give a proof for the sake of completeness.
\begin{Lemma}\label{weakly mixing}
	Let $\Gamma$ be an i.c.c.\@ group. Then the conjugation action on $L\Gamma$ is weakly mixing.
\end{Lemma}
\begin{proof}
	Let $\Gamma=\{\gamma_j|\,j\in\NN\}$ be an enumeration of $\Gamma$. Assume $ H\subset L\Gamma\subset \ell^2(\Gamma)$ is an $\Gamma$-invariant, finite-dimensional subspace and let $\{\xi_1,\dots,\xi_n\}$ be an orthonormal basis of $H$ such that $\xi_1\notin\CC\delta_e$. Then for every $\eps>0$ there is a $K\in\NN$ such that 
	\begin{align}\label{Def K}
	\left\| \xi_j-\sum_{i=1}^K\left\langle \xi_j,\delta_{\gamma_i}\right\rangle \delta_{\gamma_i}\right\| <\eps\quad\text{ for all } j=1,\dots,n.
	\end{align}
	Let $F=\{\gamma_1,\dots,\gamma_K\}$. Then by \cite[Proposition 3.4]{CTU2013} there exists a $\gamma\in\Gamma$ such that 
	\begin{align}\label{ICU}
	\gamma F\gamma^{-1}\cap F\subset\{e\}.
	\end{align} 
	Let $H_F:=\spa(F)$ and $P_F$ the orthogonal projection on $H_F$.
	
	As $\{\gamma\xi_1\gamma^{-1},\dots,\gamma\xi_n\gamma^{-1}\}$ is again an orthonormal basis of $H$ we have $c_j\in\CC$ with $\sum_{j=1}^n|c_j|^2=1$ such that
	\[\xi_1=\sum_{j=1}^nc_j \gamma\xi_j\gamma^{-1}
	=\sum_{j=1}^nc_j \left( \sum_{i=1}^K\left\langle \xi_j,\delta_{\gamma_i}\right\rangle \delta_{\gamma\gamma_i\gamma^{-1}}
	+\sum_{i=K+1}^\infty\left\langle \xi_j,\delta_{\gamma_i}\right\rangle \delta_{\gamma\gamma_i\gamma^{-1}}\right)  .	\]
	We have  $\sum_{i=1}^K\left\langle \xi_j,\delta_{\gamma_i}\right\rangle \delta_{\gamma\gamma_i\gamma^{-1}}\in H_F^\perp+\CC\delta_e$ because of (\ref{ICU}), which together with (\ref{Def K}) implies
	\begin{align*}
	\|P_F(\xi_1)\|&\leq |\left\langle \xi_1,\delta_e\right\rangle| 
	+\left\| P_F\left( \sum_{j=1}^n c_j\sum_{i=K+1}^\infty\left\langle \xi_j,\delta_{\gamma_i}\right\rangle \delta_{\gamma\gamma_i\gamma{-1}}\right) \right\|  \\
	&\leq |\left\langle \xi_1,\delta_e\right\rangle| +\eps \sum_{j=1}^n |c_j|\\
	&\leq |\left\langle \xi_1,\delta_e\right\rangle| +n\eps.	\end{align*}
	Since $\|P_F(\xi_1)\|>1-\eps$ by (\ref{Def K}), we get a contradiction when choosing $\eps<n^{-1}(1-|\left\langle \xi_1,\delta_e\right\rangle|)$.
\end{proof}
\begin{Def}
	We call an extremal character \emph{conjugation weakly mixing} if the conjugation action on its GNS construction is weakly mixing.
\end{Def}
\begin{Qu}
	Which other characters are conjugation weakly mixing?
\end{Qu}
The following statement contains Theorem \ref{disintegration rigidity} as a special case.
\begin{Thm}
	Let $\tau$ be a conjugation weakly mixing character on $\Gamma$. Then $(\Gamma,\tau)$ is disintegration rigid.
\end{Thm}
\begin{proof}
	We first assume that $\alpha$ is ergodic.
	An action on a finite von Neumann algebra $\sigma\colon \Gamma\curvearrowright N$ is weakly mixing if and only if for every action $\alpha\colon \Gamma\curvearrowright A$ on a finite von Neumann algebra one has $(A\overline{\otimes}N)^{(\alpha\otimes\sigma)}=A^\alpha\otimes 1$  \cite[Proposition D.2]{Vaes2006}. So if we take $A=L^\infty(\Omega)$ as in Section \ref{general irpdf} and $N=\pi(\Gamma)''$, Lemma \ref{weakly mixing} implies that \[M^\theta=(A\overline{\otimes}N)^{(\alpha\otimes \conj(\pi))}=A^\alpha=\CC.\]	
	
	$\tau$ is extremal because if the conjugation action is weakly mixing, it must be ergodic, hence the GNS construction is a factor.
	Hence by Proposition \ref{M theta} every i.r.p.d.f.~$\varphi$  with $\mb E[\varphi]=\tau$ is given by an element in $M^\theta$, which	proves the statement in the ergodic case.
	
	The general case follows by ergodic decomposition:  Let $\varphi$ be an i.r.p.d.f.~with $\mb E[\varphi]=\tau$. Then the restriction to the ergodic components are ergodic i.r.p.d.f.'s. The expectation values of these ergodic i.r.p.d.f.'s integrate to $\tau$ and are therefore by extremality $\mu$-almost surely equal to $\tau$. Hence we can apply the statement to them and get that they are equal to $\tau$ $\nu$-almost surely, which implies that $\varphi$ is equal $\tau$ $\mu$-almost surely.
\end{proof}
\newcommand{\etalchar}[1]{$^{#1}$}


\begin{thebibliography}{ABB{\etalchar{+}}17}

\bibitem[ABB{\etalchar{+}}11]{7s2011}
Miklos Abert, Nicolas Bergeron, Ian Biringer, Tsachik Gelander, Nikolay
  Nikolov, Jean Raimbault, and Iddo Samet.
\newblock On the growth of {B}etti numbers of locally symmetric spaces.
\newblock {\em C. R. Math. Acad. Sci. Paris}, 349(15-16):831--835, 2011.

\bibitem[ABB{\etalchar{+}}17]{7s2017}
Miklos Abert, Nicolas Bergeron, Ian Biringer, Tsachik Gelander, Nikolay
  Nikolov, Jean Raimbault, and Iddo Samet.
\newblock On the growth of {$L^2$}-invariants for sequences of lattices in
  {L}ie groups.
\newblock {\em Ann. of Math. (2)}, 185(3):711--790, 2017.

\bibitem[ADP]{ADPopa}
Claire Anantharaman-Delaroche and Sorin Popa.
\newblock {\em An introduction to {II}$_1$ factors}.
\newblock preliminary version.

\bibitem[AGN17]{AGN2017}
Miklos Abert, Tsachik Gelander, and Nikolay Nikolov.
\newblock Rank, combinatorial cost, and homology torsion growth in higher rank
  lattices.
\newblock {\em Duke Math. J.}, 166(15):2925--2964, 2017.

\bibitem[AGV14]{IRS12}
Mikl\'os Ab\'ert, Yair Glasner, and B\'alint Vir\'ag.
\newblock Kesten's theorem for invariant random subgroups.
\newblock {\em Duke Math. J.}, 163(3):465--488, 2014.

\bibitem[Bla06]{Blackadar2006}
Bruce Blackadar.
\newblock {\em Operator algebras}, volume 122 of {\em Encyclopaedia of
  Mathematical Sciences}.
\newblock Springer-Verlag, Berlin, 2006.
\newblock Theory of $C^*$-algebras and von Neumann algebras, Operator Algebras
  and Non-commutative Geometry, III.

\bibitem[CSU16]{CTU2013}
Ionut Chifan, Thomas Sinclair, and Bogdan Udrea.
\newblock Inner amenability for groups and central sequences in factors.
\newblock {\em Ergodic Theory Dynam. Systems}, 36(4):1106--1129, 2016.

\bibitem[Gel18]{Ge2018}
Tsachik Gelander.
\newblock Kazhdan-{M}argulis theorem for invariant random subgroups.
\newblock {\em Adv. Math.}, 327:47--51, 2018.

\bibitem[Hou]{Houdayer}
Cyril Houdayer.
\newblock An introduction to {II}$_1$ factors.
\newblock {\em lecture notes}.

\bibitem[L{\"u}c02]{Luck}
Wolfgang L{\"u}ck.
\newblock {\em {$L^2$}-invariants: theory and applications to geometry and
  {$K$}-theory}, volume~44 of {\em Ergebnisse der Mathematik und ihrer
  Grenzgebiete. 3. Folge.}
\newblock Springer-Verlag, Berlin, 2002.

\bibitem[Str81]{Stratila1981}
{\c S}erban Str{\u a}til{\u a}.
\newblock {\em Modular theory in operator algebras}.
\newblock Editura Academiei Republicii Socialiste Rom\^ania, Bucharest; Abacus
  Press, Tunbridge Wells, 1981.

\bibitem[Tak02]{Takesaki2002}
Masamichi Takesaki.
\newblock {\em Theory of operator algebras. {I}}, volume 124 of {\em
  Encyclopaedia of Mathematical Sciences}.
\newblock Springer-Verlag, Berlin, 2002.
\newblock Operator Algebras and Non-commutative Geometry, 5.

\bibitem[Tak03]{Takesaki}
Masamichi Takesaki.
\newblock {\em Theory of operator algebras. {II}}, volume 125 of {\em
  Encyclopaedia of Mathematical Sciences}.
\newblock Springer-Verlag, Berlin, 2003.
\newblock Operator Algebras and Non-commutative Geometry, 6.

\bibitem[Tho64a]{Thoma1964(3)}
Elmar Thoma.
\newblock {D}ie unzerlegbaren, positiv-definiten {K}lassenfunktionen der
  abz\"ahlbar unendlichen, symmetrischen {G}ruppe.
\newblock {\em Math. Z.}, 85:40--61, 1964.

\bibitem[Tho64b]{Thoma1964(1)}
Elmar Thoma.
\newblock {{\"U}}ber unit\"are {D}arstellungen abz\"ahlbarer, diskreter
  {G}ruppen.
\newblock {\em Math. Ann.}, 153:111--138, 1964.

\bibitem[Vae07]{Vaes2006}
Stefaan Vaes.
\newblock Rigidity results for {B}ernoulli actions and their von {N}eumann
  algebras (after {S}orin {P}opa).
\newblock {\em Ast\'erisque}, (311):Exp. No. 961, viii, 237--294, 2007.
\newblock S\'eminaire Bourbaki. Vol. 2005/2006.

\bibitem[VK81]{VerKer1981}
Anatoly Vershik and Sergei Kerov.
\newblock Characters and factor representations of the infinite symmetric
  group.
\newblock {\em Dokl. Akad. Nauk SSSR}, 257(5):1037--1040, 1981.

\end{thebibliography}
\end{document}